\documentclass[12pt,reqno]{amsart}
\usepackage{amssymb}
\usepackage{amsmath}
 \usepackage[usenames,dvipsnames]{color}
\usepackage[kerning=true]{microtype}
\usepackage[page]{appendix}
 \usepackage{pgf,tikz}
 \usetikzlibrary{decorations.markings}
\usetikzlibrary{arrows,shapes,positioning,calc,shadows}
\tikzstyle arrowstyle=[scale=1]
\tikzstyle directed=[postaction={decorate,decoration={markings,
    mark=at position .5 with {\arrow[arrowstyle]{stealth}}}}]

\usepackage{fullpage}
\usepackage{amsthm}
\usepackage{bm}
\usepackage{shuffle}
\usepackage[T1]{fontenc}
\usepackage{graphicx}
\usepackage[applemac]{inputenc}
\usepackage[french]{todonotes}

\usepackage{wrapfig}

\usepackage[colorlinks=true,
citecolor=cyan,urlcolor=blue,linktocpage=true,pagebackref]{hyperref}

\def\bleu{\textcolor{blue}}

\def\auteur#1{{\sc #1}}
\def\titreref#1{{\em #1}}
\def\vol#1{{\bf #1}}

\newtheorem{conjecture}{\bleu{Conjecture}}

\newtheorem{prop}{\bleu{Proposition}}
\newtheorem{cor}{\bleu{Corollary}}

\numberwithin{equation}{section}
\numberwithin{lemma}{section}
\numberwithin{rmk}{section}
\numberwithin{cor}{section}

\newcommand{\define}[1]{\bleu{\bf{#1}}}
\def\ie{{\it i.e.}~}

\def\MR#1{\href{http://www.ams.org/mathscinet-getitem?mr=#1}{MR#1}}

\def\pref#1{{\rm (\ref{#1})}}

\newcommand{\BF}{\mathbb{DBF}}
\newcommand{\Bos}{\bm{x}}
\newcommand{\qbinom}[2]{\genfrac{[}{]}{0pt}{}{#1}{#2}_q}

\newcommand{\Fer}{\bm{\theta}}
\newcommand{\GL}{\mathrm{GL}}

\newcommand{\stir}[2]{\genfrac{\{}{\}}{0pt}{}{#1}{#2}}

\renewcommand{\S}{\mathbb{S}}

\newcommand{\Q}{\mathbb{Q}}
\newcommand{\Real}{\mathbb{R}}
\newcommand{\R}{\mathcal{R}}
\newcommand{\G}{\mathcal{G}}

\title[Bosons Fermions]{The Bosonic-Fermionic Diagonal Coinvariant Modules Conjecture
}
\author{F.~Bergeron}
\address{\href{http://bergeron.math.uqam.ca}{D\'epartement de Math\'ematiques, Lacim, UQAM.}}
  \email{\href{mailto:bergeron.francois@uqam.ca}{bergeron.francois@uqam.ca}}
  \date{\bleu{\bf \today}. This work was supported by NSERC}

\begin{document}


\begin{abstract} We describe a general conjecture on how one may derive from the generic bosonic case all structural properties of multivariate diagonal coinvariant modules in any number (say $k$) of sets of $n$ commuting variables (bosons), and any number (say $j$) of sets of $n$ anticommuting variables (fermions).   \end{abstract}

\maketitle
 \parskip=0pt
\footskip=30pt

\parskip=8pt  
\parindent=20pt

\section{Introduction}

Much interesting work has been done recently on diagonal coinvariant spaces spaces in both commuting and anticommuting variables. See for instance~\cite{orellana,swanson,swanson_wallach,wallach}. The purpose of this short note is to present a general conjecture expressing the fact that one can simply calculate all cases of multivariate diagonal coinvariant modules in $k$ sets of $n$ commuting variables (bosons), and $j$ sets of $n$ anticommuting variables (fermions), just from the generic case of multivariate diagonal coinvariant spaces.   

\section{\bleu{Global setup}}
Let $\Bos$ and $\Fer$ respectively be matrices  of variables; with $\Bos=(x_{ab})$ a $k\times n$ matrix, and $\Fer=(\theta_{cd})$ a $j\times n$ matrix. One may even assume that $k$ and $j$ are infinite.
The variables in $\Bos$ commute with all variables (both in $\Bos$ and $\Fer$), whereas the variables in $\Fer$ are anticommuting among themselves, \ie for $\theta$ and $\theta'$ in $\Fer$ one has $\theta\theta'=-\theta'\theta$. The $\Bos$-variables are said to be \define{bosomic}, and those in $\Fer$ to be \define {fermionic}.
We consider that the ring of polynomials $\R=\R_n=\R_{k,j;n}:=\Real[\Bos;\Fer]$ comes equipped with the group action (expressed here with matrix multiplication)
      $$f(\Bos;\Fer) \longmapsto f(P\,\Bos\,\sigma;Q\,\Fer\,\sigma),$$
  with $P$ and $Q$ lying respectively in $\GL_k$ and $\GL_j$, and $\sigma\in\S_n$ considered here as a $n\times n$ permutation matrix. Observe that $\sigma$ acts similarly on both the $\Bos$ and $\Fer$ variables, by permuting columns. It has become usual to say that this is a diagonal action of $\S_n$. It is worth underlining that the three single actions commute.
 
 As is often done, we denote by $\R^{\S_n}$ the subring of \define{invariants} of $\R$, \ie the polynomials that are invariant under the (single) action of $\S_n$. The associated \define{coinvariant} module (often denoted by $\R_{\S_n}$, but here denoted otherwise for reasons that will become clear in the sequel) is then set to be the quotient
       $$\BF_{k,j;n}:=\R_{k,j;n}/\langle \R^{\S_n}_+\rangle, $$
    with $\R^{\S_n}_+$ standing for the constant term free portion of $\R^{\S_n}$. Since $\R^{\S_n}$ is globally invariant under the action of $\G=\GL_k\times \GL_j\times \S_n$, there is an induced action of $\G$ on $\BF_{k,j;n}$.
  
Consider any $\mathcal V$, which is a $\G$-submodule (or stable quotient module)  of $\R=\R_{k,j;n}$.  As is well known, the decomposition  of $\mathcal{V}$ into irreducibles of  is entirely encoded in the symmetric function expression\footnote{Observe that various sets of variables are separated by semi-colons.}
\begin{equation}\label{charactere}
   \mathcal V(\bm{q};\bm{u};\bm{z}):=\sum_{\mu\vdash n} \sum_{\lambda,\rho} v_{\lambda\rho\mu} s_\lambda(\bm{q}) s_\rho(\bm{u}) s_\mu(\bm{z}),
 \end{equation}
where  $s_\lambda(\bm{q})$ and  $s_\rho(\bm{u})$ are respectively characters for (polynomial) irreducible representations of $\GL_k$ and $\GL_j$ expressed as functions of $\bm{q}=q_1,\ldots,q_k$ and $\bm{u}=u_1,\ldots,u_j$; and $s_\mu(\bm{z})$ is the Frobenius transform of a $\S_n$-irreducible, with $\bm{z}=z_1,z_2,\ldots$. The \define{hilbert series}\footnote{Observe that, in our notation for the passage to the Hilbert series, we simply "drop" the variables $\bm{z}$.} (or \define{graded character}) of $\mathcal{V}$ is:
\begin{equation}\label{hilbert}
   \mathcal V(\bm{q};\bm{u}):=\sum_{\mu\vdash n} \sum_{\lambda,\rho} v_{\lambda\rho\mu} s_\lambda(\bm{q}) s_\rho(\bm{u}) f^\mu,
 \end{equation}
where $f^\mu$ is the dimension of the irreducible associated to $s_\mu$, which is well known to be equal to the number of standard tableaux of shape $\mu$. Observe that the passage to the Hilbert series $\mathcal V(\bm{q};\bm{u})$ is obtained by removing $\bm{z}$ in $\mathcal V(\bm{q};\bm{u};\bm{z})$.
For the modules that we will consider, the coefficients $ c_{\lambda\rho\mu}$ do not depend on $k$ and $j$. The dependence on $k$ and $j$ is rather reflected in the fact that some of the functions $s_\lambda(\bm{q})$ and  $s_\rho(\bm{u})$ when the number of variables is too small, {\sl i.e.} $k$ (resp. $j$) is less than the number of parts of $\lambda$ (resp. $\rho$).
In other words, the stable expression for $ \mathcal V(\bm{q};\bm{u};\bm{z})$ is obtained whenever $k$ and $j$ become large enough\footnote{It is sufficient to take $k$  larger or equal to $n$, since this holds for the whole space of polynomials.}. Such modules are said to be \define{coefficient stable}. 

When this is the case, it is often useful to write \pref{charactere} in the form of a ``variable free'' expression:
   $$   \mathcal V:=\sum_{\mu\vdash n} \sum_{\lambda,\rho} v_{\lambda\rho\mu}\, s_\lambda\otimes s_\rho\otimes s_\mu.$$
Using plethystic notation\footnote{See~\cite{bergeron} for notions not described here.}, 
\begin{prop} One has
\begin{align}
   \R_n(\bm{q};\bm{u};\bm{z})&=h_n[\Omega[\bm{q}-\varepsilon\,\bm{u}]\,\bm{z}]\\
       &=\sum_{\mu\vdash n} s_\mu[\Omega[\bm{q}-\varepsilon\,\bm{u}]]\, s_\mu(\bm{z}),
\end{align}
   with $\Omega=\sum_{i\geq 0} h_i$. 
\end{prop}
For the above plethystic calculation, $\varepsilon$ is defined to be such that $p_i[\varepsilon]=(-1)^{i}$, so that $p_i[-\varepsilon\,\, \bm{u}]=\omega\, p_i(\bm{u})$. In particular, by the classical summation formula for Schur functions, one gets
   $$s_\theta[\bm{q}-\varepsilon\,\bm{u}] = \sum_{\nu\subseteq \theta} s_\nu(\bm{q}) s_{\theta'/\nu'}(\bm{u}).$$
Exploiting the Hall scalar product on symmetric function\footnote{For which the Schur functions $s_\mu$ form an orthonormal basis.} in $\bm{z}$, 
we may consider the \define{coefficient} $R_\mu(\bm{q};\bm{u}):=\langle \R_n(\bm{q};\bm{u};\bm{z}),s_\mu(\bm{z})\rangle$ of $s_\mu$ in $\R_n$. As a variable free expression, we thus set
    $$R_\mu:=\langle \R_n,s_\mu\rangle=\sum_{\lambda,\rho} d_{\lambda\rho\mu}\, s_\lambda\otimes s_\rho,$$
   and likewise for any coefficient stable $\G$-module $\mathcal{V}$:
       $$V_\mu:=\langle \mathcal{V}_n,s_\mu\rangle=\sum_{\lambda,\rho} v_{\lambda\rho\mu}\, s_\lambda\otimes s_\rho,$$
  It is clear that the $R_\mu$ form an upper bound for all $V_\mu$, so that 
     $$0\leq v_{\lambda\rho\mu}\leq d_{\lambda\rho\mu}, \qquad \hbox{for all } \lambda,\rho,\ {\rm and}\ \mu.$$
  Thus it is interesting to observe that the Cauchy kernel formula implies that   
 \begin{cor}
   When $j=0$, the coefficient $R_\mu(\bm{q})=\sum_{\lambda} d_{\lambda\mu} s_\lambda(\bm{q})$ of $s_\mu(\bm{z})$ in $\R_n(\bm{q};\bm{z})$ is given by the formula
     \begin{align}
        R_\mu(\bm{q}) &= s_\mu[\Omega(\bm{q})] \nonumber\\
                   &= s_\mu[\textstyle \sum_{i\geq 0} h_i(\bm{q})].
        \end{align}
  \end{cor}
 
 \section{\bleu{The boson-fermion modules}}
 As show in \cite{multicomplex}, there exist a coefficient stable expression for the pure bosonic (commuting variables) multivariate coinvariant module, which we denote by
 $$ \mathcal{E}_n=\sum_{\mu\vdash n}\mathcal{C}_\mu\otimes s_\mu,\qquad {\rm with}\qquad
    \mathcal{C}_\mu:=\sum_\lambda c_{\lambda\mu}\,s_\lambda.$$
The integers $c_{\lambda\mu}$ are non-vanishing only for partitions $\lambda$ of size at most $\binom{n}{2}-\eta(\mu')$, and having at most $n-\mu_1$ parts.  Recall that $\eta(\mu):=\sum_{i} \mu_i\,(i-1)$.
Expressed in terms of variables, the above expression takes the form
\begin{equation}\label{bosons}
   \mathcal{E}_n(\bm{q};\bm{z})= \sum_{\mu\vdash n}\sum_\lambda c_{\lambda\mu}\,s_\lambda(\bm{q}) s_\mu(\bm{z}).
\end{equation}
Our main conjecture is that
\begin{conjecture} [Diagonal Supersymmetry]\label{mainconj} The multigraded Frobenius characteristic of the boson-fermion diagonal modules may be calculated from the generic Frobenius characteristic for bosons modules via the universal formula
 \begin{align}\label{boson_fermion_equation}
      \BF_{k,j;n}(\bm{q};\bm{u};\bm{z}) &=  \mathcal{E}_n(\bm{q}-\varepsilon\,\bm{u};\bm{z})\\
      &= \sum_{\mu\vdash n} \mathcal{C}_{\mu}[\bm{q}-\varepsilon\, \bm{u}] s_\mu(\bm{z})\\
      &=\sum_{\mu\vdash n}\sum_\lambda c_{\lambda\mu}\,s_\lambda[\bm{q}-\varepsilon\,\bm{u}] s_\mu(\bm{z}).
   \end{align}
\end{conjecture}
\noindent
Thus, the $(k,j)$-multi-degree enumeration (or $\GL_k\times \GL_j$-character) of the $\S_n$-irreducible component of type $\mu$ in  $\BF_{n}$ is obtained as
 \begin{align}\label{boson_fermion_equation_coeff}
      \mathcal{C}_{\mu}[\bm{q}-\varepsilon \bm{u}]  
      &=\sum_\lambda c_{\lambda\mu}\,s_\lambda[\bm{q}-\varepsilon\,\bm{u}]\\
      &= \sum_\lambda c_{\lambda\mu} \sum_{\nu\subseteq \lambda} s_\nu(q_1,\ldots,q_k) 
          s_{\lambda'/\nu'}(u_1,\ldots,u_j).
  \end{align}
Observe that the specification of $k$ and $j$ in $\BF_{k,j;n}(\bm{q};\bm{u};\bm{z})$  is redundant once the parameters $\bm{q}=q_1,\ldots,q_k$ and $\bm{u}=u_1,\ldots,u_j$ are specified. We may thus omit them when this is the case, and also formaly write 
\begin{equation}\label{KkJj}
	\BF_n=\sum_{\mu\vdash n}\sum_\lambda \sum_{\nu\subseteq \lambda}c_{\lambda\mu} \, s_\nu\otimes 
s_{\lambda'/\nu'}\otimes s_\mu,
\end{equation} 
for the \define{generic diagonal Boson-Fermion Frobenius}. 
Clearly, we have 
\begin{align}
      &\BF_n(q;0;\bm{z})=h_n^*(\bm{z})/h_n^*(1),\qquad {\rm with}\qquad f^*(\bm{z}):=f[\bm{z}/(1-q)],\label{K1J0}\\[4pt]
      &\BF_n(0;u;\bm{z})=\sum_{a=0}^{n-1} u^a s_{(n-a,1^a)}(\bm{z}),\label{K0J1}\\[4pt]
      &\BF_n(q,t;0;\bm{z})=\nabla(e_n)(q,t;\bm{z}),\label{K2J0}\\[4pt]     
      &\BF_n(q_1,\ldots,q_k;0;\bm{z})=\mathcal{E}_n(q_1,\ldots,q_k\bm{z}),\label{KkJ0}\\[4pt]
      &\BF_n(0;u_1,\ldots,u_j;\bm{z}) = \sum_{\mu\vdash n} \sum_\lambda c_{\lambda\mu}\,s_{\lambda'}(u_1,\ldots,u_j) s_\mu(\bm{z}) \label{K0Jj}
 \end{align}
The $\nabla$ operator occurring in~\pref{K2J0} is a Macdonald ``eigenoperator'' (introduced in~\cite{SF}). This is to say it affords as eigenfunctions the (modified) Macdonald operators, usually denoted $\widetilde{H}_\mu$.

 It has been conjectured\footnote{Notice that one of the parameters is equal to $1$. This is because the lacking ``statistic'' on Dyck-path pairs $(\alpha,\beta)$ is not yet known.} in~\cite{trivariate} that
 \begin{equation}\label{K3J0}
   \BF_n(q,t,1;0;\bm{z})=\sum_{\alpha\preceq\beta} q^{\rm{dist}(\alpha,\beta)}\,\mathbb{L}_\beta(t;\bm{z}),
 \end{equation}
 where the sum is over pairs of elements of the Tamari lattice, and $\rm{dist}(\alpha,\beta)$ is the length of the longest chain going from $\alpha$ to $\beta$. Here, $\mathbb{L}_\beta(t;\bm{z})$ stands for the LLT-polynomial associated to the Dyck-path $\beta$ (see~\cite{GLk_x_Sn} for more details). Furthermore, it has been conjectured by N.~Bergeron-Machacek-Zabrocki\footnote{In fact, this follows from~\pref{K2J1}, via a formula of Haglund, Rhoades and Shimozono (see~\cite{HRS}).} that
\begin{equation}\label{K1J1}
    \BF_n(q;u;\bm{z}) = \sum_{k=0}^{n-1} \sum_{\lambda\vdash n}\sum_{\tau\in {\rm SYT}(\lambda)} q^{\alpha(\tau)}\qbinom{{\rm des}(\tau)}{k}\, u^{k}\,s_{\lambda}(\bm{z}),
\end{equation}
where, for a standard tableau $\tau$ of shape $\lambda$, one sets 
  $$\alpha(\tau) := {\rm maj}(\tau)-k\,{\rm des}(\tau)+\binom{k}{2}.$$
In~\cite{rhoades}, Kim and Rhoades show that
 \begin{equation}\label{K0J2}
    \BF_n(0;u,v;\bm{z}) = \sum_{a+b\,\leq\, n-1} u^a\,v^b \, \big(s_{(n-a,1^a)} \star s_{(n-b,1^b)}
        -s_{(n-(a-1),1^{a-1})} \star s_{(n-(b-1),1^{b-1})}\big),
\end{equation}
   with ``$\star$'' standing for the Kronecker product. Denoting by $g_{\alpha,\beta}^\mu$ the  \define{Kronecker coefficients}:
      $$g_{\alpha,\beta}^\mu:=\langle s_\alpha\star s_\beta,s_\mu\rangle,$$
   one may reformulate the above as
  \begin{equation}
    \BF_n(0;u,v;\bm{z}) = \sum_{\mu\vdash n}  \left(\sum_{b+d\,\leq\, n-1} u^b\,v^d \, ( g_{(a\,|\,b),(c\,|\,d)}^\mu-g_{(a+1\,|\,b-1),(c+1\,|\,d-1)}^\mu)  \right)\,s_\mu(\bm{z}),
\end{equation}
using the Frobenius notation $(a\,|\,b) = (a+1,1^b)$ for hook-shaped partitions. For each term in the inner sum above, we assume that $a+b=n-1$ (likewise for $c$ and $d$). The differences are know to be positive (see~\cite{remmel}). The various results (see~\cite{snowden}) on the stability of Kronecker coefficients certainly have a bearing here, since they imply corresponding stabilities for the coefficients of the $s_\mu$.  
   
Zabrocki has conjectured (see~\cite{zabrocki}) that 
   \begin{equation}\label{K2J1}
      \BF_n(q,t;u;\bm{z}) =\sum_{a=0}^{n-1} u^{a} \Delta'_{e_{n-a-1}}(e_n(\bm{z})).
   \end{equation}
   The parameters $q$ and $t$ arise from the application of the operators $\Delta'_{e_k}$.
 It follows that \pref{K1J1} may also be written as
 \begin{equation}
    \BF_n(q;u;\bm{z}) = \sum_{a=0}^{n-1} u^{a} \Delta'_{e_{n-a-1}}(e_n(\bm{z}))\Big|_{t=0}.
\end{equation}
Finally, we have
\begin{equation}\label{K1J2}
    \BF_n(1;2;\bm{z}) = \frac{1}{2}\sum_{\mu\vdash n} 2^{\ell(\mu)} (-1)^{n-\ell(\mu)} \binom{\ell(\mu)}{d_1,\ldots,d_n} \,p_\mu(\bm{z}),
\end{equation}
where $d_i=d_i(\mu)$ stands for the number of parts of of size $i$ in $\mu$. Finally, D'Adderio-Iraci-Wyngaerd conjecture in \cite[Conj. 8.2.]{dadderio} the more inclusive identity:
\begin{equation}\label{K2J2}
    \BF_n(q,t;u,v;\bm{z}) = \sum_{k=1}^{n-1}\sum_{i+j =k} u^i v^j \Theta_{e_ie_j} \nabla(e_{n-k}),
\end{equation}
where, for any symmetric functions $g$ and $f$, $\Theta_g f$ is defined as
     $$\Theta_g f(\bm{z}) := \Pi\, g^*\,\Pi^{-1} f(\bm{z}),\qquad
        {\rm setting}\qquad g^*(\bm{z}):=g[\bm{z}/(1-t)(1-q)] .$$
Here, $\Pi$ stands for the Macdonald eigenoperator having as eigenvalues for $\widetilde{H}_\mu$ the product $\prod_{(i,j)\in\mu/(1)} (1 - q^it^j)$, for $(i,j)$ running over cartesian coordinates of cells in $\mu$ (omitting the cell $(0,0)$).

	$$\mathcal{G}(q,t;\underline{z};\underline{x}) := \sum_{k=0}^{n-1} \sum_{\nu\vdash k} m_\nu(\underline{z})\, \Theta_{e_\nu} \nabla(e_{n-k})(q,t;\underline{x}).$$

Table~\ref{table0}  summarizes the overall situation\footnote{With $k$ standing for the numbers of sets of commuting variables, and $j$ for those that are anticommuting.}.
\begin{table}[ht]
\begin{center}
\begin{tabular}{|c||c|c|c|c|c|}
\hline
$k\setminus j$ & 0 &  1 & 2 & $\cdots$ & $j$ \\
\hline
\hline
0  &  1 & \pref{K0J1}& \pref{K0J2} &$\cdots$&\pref{K0Jj}\\
\hline
 1  & \pref{K1J0}  & \pref{K1J1} &\pref{K1J2} &$\cdots$&\\
 \hline
 2  & \pref{K2J0} & \pref{K2J1} &\pref{K2J2} &$\cdots$& \\
\hline
3  & \pref{K3J0} & & &$\cdots$& \\
\hline
$\vdots$  &\vdots  & \vdots&\vdots&$\ddots$&\vdots \\
\hline
$k$  & \pref{KkJ0} & & &$\cdots$&\pref{KkJj}\\
\hline
\end{tabular}\vskip5pt\caption{Overall situation for the various formulas.}\label{table0}
\end{center}
\vskip-10pt
\end{table}
Conjecture~\ref{mainconj} essentially states all entries may be obtained from \pref{KkJ0} (or equivalently from \pref{K0Jj}).
It is interesting to observe that one obtains polynomial expressions in $k$ and $j$, when setting all parameters $q_i=1$ and $u_j=1$. More precisely, writing  $\BF_n(k;j;\bm{z})$ 
for the resulting expression,
we have the following.
\begin{prop}
The coefficients of each $s_\mu(\bm{z})$, in the Schur expansion of $\BF_n(k;j;\bm{z})$, is a polynomial in $k$ and $j$, with coefficients in $\Q$. Hence, this is also the case for the associated dimension\footnote{Obtained by replacing each $s_\mu(\bm{z})$ by the number, $f^\mu$, of standard tableaux of shape $\mu$.} $\BF_n(k;j)$.
\end{prop}

\section{Links with the main conjecture}
Conjecture~\pref{K2J1} directly led to our main conjecture, in view of an elegant link (first stated in 2017, but only recently published) between the generic expression for $\mathcal{E}_n$ and the effect of the $\Delta'_{e_k}$ operators on $e_n$. The precise relevant statement (see~\cite[Conj. 1]{GLk_x_Sn}) says that
\begin{conjecture}[Delta via skew]\label{conj_skew}
For all $k$,
\begin{align}
   (e_k\otimes \mathrm{Id})\,\mathcal{E}_n&= \sum_{\mu\vdash n} (e_k^\perp \mathcal{C}_{\mu})(q,t)\, s_\mu(\bm{z}) \nonumber\\ 
                  & = \Delta'_{e_{n-k-1}}(e_n(\bm{z})).
 \end{align}
 \end{conjecture}
\noindent In other words, we get  $\Delta'_{e_{n-k-1}}(e_n(\bm{z}))$ from $\mathcal{E}_n$, first by applying the skew operator $e^\perp_{k}$ to the various $s_\lambda$, and then by evaluation of the resulting expression in $q$, $t$. To see how this relates to our general conjecture, we recall that effect on a symmetric function $f(\bm{q})=f(q_1,\ldots,q_k)$ of the operator $\sum_{a=0}^{n-1} u^a e^\perp_k$ may be globally expressed in plethystic notation as
     $$\sum_{a=0}^{n-1} u^a e^\perp_k f(\bm{q}) = f[\bm{q}-\varepsilon\, u].$$
 Thus, assuming that Conjecture~\ref{conj_skew} holds, we see that~\pref{K2J1} may be now simply be coined as
    $$\BF_n(q,t;u;\bm{z}) = \mathcal{E}_n[q+t-\varepsilon\, u;\bm{z}].$$
 This immediately\footnote{During the January 2019 Banff meeting where Mike Zabrocki presented his conjecture for the first time.} suggested that the more general formula of Conjecture~\ref{mainconj} should hold. All experiments confirmed this. 
 Known or conjectured formulas (due to various researchers) for the dimensions of $\BF_n(k;j)$ as functions of $n$, for small $k$ and $j$, are displayed in Table~\ref{table1}.
 \begin{table}[h!]
\begin{center}
\def\arraystretch{1.5}
\begin{tabular}{|c||c|c|c|c|c|}
\hline
$k\setminus j$ & 0 &  1 & 2 \\
\hline
\hline
0  &  1 & $2^{n-1}$ &$\binom{2n-1}{n}$ \\[4pt]
\hline
 1  & $n!$  & $\sum_{i=1}^n i! {\textstyle\stir{n}{i}}$ & $2^{n-1}\,n!$ \\
 \hline
 2  & $(n+1)^{n-1}$ & $\sum_{i=0}^{n+1} \binom{n+1}{i}\frac{i^n}{2(n+1)}$ & ?  \\
  \hline
3  & $2^n(n+1)^{n-2}$ & ? & ? \\
\hline
\end{tabular}\vskip5pt\caption{Dimensions of $\BF_n(k,j)$.}\label{table1}
\end{center}\vskip-15pt
\end{table}
Here, $\stir{n}{k}$ stands for the Stirling numbers of the second kind.  Currently known or conjectured formulas for the multiplicities of alternating component in $\BF_n$ appear in Table~\ref{table2},
 \begin{table}[h!]
\begin{center}
\def\arraystretch{1.5}
\begin{tabular}{|c||c|c|c|c|c|}
\hline
$k\setminus j$ & 0 &  1 & 2  & 3\\
\hline
\hline
0  &  0 & $1$ &$n$  & $n^2-n+1$\\[4pt]
\hline
 1  & $1$  & $2^{n-1}$ & $3^{n-1}$ & $2^{-1} F_{3n-1}$\\
 \hline
 2  & $\frac{1}{n+1}\binom{2n}{n}$ & $s(n)$ & $\frac{2^{n-1}}{n+1}\binom{2n}{n}$ & ? \\
  \hline
3  & $\frac{2}{n(n+1)} \binom{4\,n+1}{n-1}$ & ? & ?  & ?\\
\hline
\end{tabular}\vskip5pt\caption{Coefficients of $s_{1^n}(\bm{z})$ in $\BF_n(k,j)(\bm{z})$.}\label{table2}
\end{center}\vskip-15pt
\end{table}
where $s(n)=\frac{1}{n} \sum_{i=0}^{n-1} \binom{n}{i}\, \binom{n}{i+1}\,2^i$ denotes the $n^{\rm th}$ small Schröder number, and $F_n$ stands for the $n^{\rm th}$ Fibonacci number.

Formulas for low degree components of the corresponding general expressions have also been conjectured to hold. The first of these, see~\cite{multicomplex}, states that
    \begin{equation}\label{Conj_low_degree}
    	\BF_n(\bm{q};0;\bm{z})=_{(n)} \frac{h_n[\Omega(\bm{q})\,\bm{z}]}{h_n[\Omega(\bm{q})]} ,
    \end{equation}
 where $\Omega(\bm{q})=1+h_1(\bm{q})+h_2(\bm{q})+\ldots $; with equality holding for terms of degree at most $n$ in the $\bm{q}$ variables. It may be worth recalling here that, when $\bm{q}$ consists of only one variable, the above is well known to be an equality. Indeed, this the symmetric group case of the Chevalley-Shephard-Todd theorem. 
 
A slightly stronger form of a conjecture stated in~\cite{dadderio}, is that the difference $\BF_n(q,t;\bm{u};\bm{z}) -\mathcal{M}_n(q,t;\bm{u};\bm{z}) $ is Schur positive in all three sets of variables $\bm{u}=\{u_1,u_2,\ldots,u_r\}$ (for all $r$), $\bm{z}$, and $\{q,t\}$; setting
     \begin{equation}\label{Conj_low_theta}
    	\mathcal{M}_n(q,t;\bm{u};\bm{z}) := \sum_{k=0}^{n-1} \sum_{\nu\vdash k} m_\nu(\bm{u})\, \Theta_{e_\nu} \!\nabla(e_{n-k})(q,t;\bm{z}).
    \end{equation}
Furthermore, the expression $\mathcal{M}_n(q,t;\bm{u};\bm{z})$ is itself Schur positive in all three sets of variables.

Although some of the terms in either of  the expressions (right hand-side of, \ref{Conj_low_degree}) and \pref{Conj_low_theta} may be recovered (via plethysm) from the other, this not the case for all terms. Hence, the two statements are not equivalent.

\section{An explicit example}
With $n=3$, we have 
$$\mathcal{E}_3=1\otimes s_{3}+(s_{1} + s_{2})\otimes s_{21}+(s_{11} + s_{3})\otimes s_{111},$$
from which we deduce that
  $$\begin{aligned}
 \BF_3&=1 \otimes 1 \otimes s_{3} 
 + ( s_{1} \otimes 1  + 1 \otimes s_{1} +  s_{2} \otimes 1 + s_{1} \otimes s_{1}  + 1 \otimes s_{11}     ) \otimes s_{21} \\
&+ (s_{11} \otimes 1 + s_{1} \otimes s_{1}   + 1 \otimes s_{2} + s_{3} \otimes 1 +\  s_{2} \otimes s_{1}+ s_{1} \otimes s_{11}   + 1 \otimes s_{111}     ) \otimes s_{111}.
             \end{aligned}$$
By specialization, we get
\begin{align*}
      &\BF_3(q;0;\bm{z})= s_{3}(\bm{z})+(q+q^2)\, s_{21}(\bm{z})+q^3\, s_{111}(\bm{z}), \\[4pt]
      &\BF_3(0;u;\bm{z})=  s_{3}(\bm{z}) + u\,s_{21}(\bm{z})+ u^2\, s_{111}(\bm{z}), \\[4pt]
      &\BF_3(q,t;0;\bm{z})=s_{3}(\bm{z})+(q^{2} + q t + t^{2} + q + t)\,s_{21}(\bm{z}) \\
         &\qquad\qquad\qquad\qquad\qquad +(q^{3} + q^{2} t + q t^{2} + t^{3} + q t)\,s_{111}(\bm{z}),  \\[4pt]     
      &\BF_3(q,t;u;\bm{z})=  s_{3}(\bm{z})  + ( q+t  + u +  q^2+q\,t+t^2 + q\,u+t\,u ) \, s_{21}(\bm{z}) \\
			&\qquad + (q\,t + q\,u+t\,u + u^2 + q^3+q^2t+q\,t^2+t^3+q^{2}u + q\, t\,u + t^{2}u) \, s_{111}(\bm{z}),\\
      &\BF_3(\bm{q};0;\bm{z})=s_{3}(\bm{z})+(s_{1}(\bm{q}) + s_{2}(\bm{q}))\,s_{21}(\bm{z})+(s_{11}(\bm{q}) + s_{3}(\bm{q}))\,s_{111}(\bm{z}), \\[4pt]
      &\BF_3(0;\bm{u};\bm{z})= s_{3}(\bm{z})
 + (s_{1}(\bm{u}) +   s_{11}(\bm{u}))\,s_{21}(\bm{z})
+ ( s_{2}(\bm{u}) +  s_{111}(\bm{u})  )\, s_{111}(\bm{z}).
 \end{align*}
 The polynomial expressions in $k$ and $j$ of the dimension, and the Frobenius characteristic for $\BF_3$ are respectively:
\begin{align*}
&\BF_3(k;j) = \frac{1}{6}  (k + j + 1) \, (k^{2} + 2\,k\, j+ j^{2} + 11\,k + 5\,j + 6),\qquad {\rm and}\\
&\BF_3(k;j;\bm{z})  =s_{3}(\bm{z})+ \frac{1}{2} (k^{2} + 2\, k\,j + t^{2} + 3\, k + j)\,s_{21}(\bm{z})\\
   &\qquad\qquad \qquad  +\frac{1}{6} (k^{3} + 3\, k^{2} j + 3\, k\,j^{2} + j^{3} + 6\,  k^{2} + 6\,  k\,j - k + 5\,  j)\,s_{111} (\bm{z}). 
\end{align*}
Explicit values for $n=3,4,$ and $5$ are the dimensions
\begin{center}
    \def\arraystretch{1.5}
    \begin{tabular}{|c||c|c|c|c|c|}
    \hline
    $k\setminus j$ & 0 &  1 & 2 \\
    \hline
    \hline
    0  &  1 & 4&10 \\
    \hline
     1  & 6  & 13 & 23\\
     \hline
     2  & 16 & 28 &  45\\
      \hline
    3  & 32 & 50&74 \\
    \hline
    \end{tabular}\qquad 
       \begin{tabular}{|c||c|c|c|c|c|}
    \hline
    $k\setminus j$ & 0 &  1 & 2 \\
    \hline
    \hline
    0  &  1 & 8&35 \\
    \hline
     1  & 24  & 75 & 192\\
     \hline
     2  & 125 & 288 &  597\\
      \hline
    3  & 400 & 785&1440 \\
    \hline
    \end{tabular}\qquad 
       \begin{tabular}{|c||c|c|c|c|c|}
    \hline
    $k\setminus j$ & 0 &  1 & 2 \\
    \hline
    \hline
    0  &  1 & 16&126 \\
    \hline
     1  & 120  & 541 & 1920\\
     \hline
     2  & 1296 & 3936 &  10541\\
      \hline
    3  & 6912 & 17072&38912 \\
    \hline
    \end{tabular}
  \end{center}

\renewcommand{\refname}{\bleu{References}}

\end{document}